\documentclass[a4paper,oneside,10pt]{amsart}
\usepackage{mathrsfs}
\usepackage{pifont}
\usepackage{amsmath}
 \usepackage{geometry}
\usepackage{amssymb}
\usepackage{graphicx}
\usepackage{geometry}
\newtheorem{theorem}{Theorem}[section]

\newtheorem{proposition}[theorem]{Proposition}

\theoremstyle{definition}
\newtheorem{definition}[theorem]{Definition}

\theoremstyle{remark}

\numberwithin{equation}{section}

\usepackage{hyperref}

\begin{document}
\begin{center}

{\huge\bf{A Note on the $_{2}F_{1}$ Hypergeometric Function}}\\[20pt]
%
{\Large Armen Bagdasaryan}
\par Institution of the Russian Academy of Sciences, 
V.A. Trapeznikov Institute for Control Sciences
\par 65 Profsoyuznaya, 117997, Moscow, Russia \par
\quad E-mail: abagdasari@hotmail.com
\\[15pt]

\end{center}



\textbf{Abstract}\par
The special case of the hypergeometric function $_{2}F_{1}$ represents the binomial series $(1+x)^{\alpha}=\sum_{n=0}^{\infty}\left(\:\begin{matrix}\alpha\\n\end{matrix}\:\right)x^{n}$
that always converges when $|x|<1$. Convergence of the series at the endpoints, $x=\pm 1$, depends on the values of $\alpha$ and needs to be checked in every
concrete case. In this note, using new approach, we reprove the convergence of the hypergeometric series for $|x|<1$ and obtain new result on its convergence at point $x=-1$ for every integer $\alpha\neq 0$, that is we prove it for the function $_{2}F_{1}(\alpha,\beta;\beta;x)$.
The proof is within a new theoretical setting based on a new method for reorganizing the integers and on the original regular method for summation of divergent series. 
\par


\par


\section{Introduction}\par

Almost all of the elementary functions of mathematics are either hypergeometric or ratios of hypergeometric functions. 
Moreover, many of the non-elementary functions that arise in mathematics and physics also have representations as hypergeometric series, that is
as special cases of a series, generalized hypergeometric function with properly chosen values of parameters
$$
_{p}F_{q}(a_{1},...,a_{p}; b_{1},...,b_{q};x)=\sum_{n=0}^{\infty}\frac{(a_{1})_{n}...(a_{p})_{n}}{(b_{1})_{n}...(b_{q})_{n}n!}x^{n}
$$
where $(w)_{n}\equiv \frac{\Gamma(w+n)}{\Gamma(w)}=w(w+1)(w+2)...(w+n-1), (w)_{0}=1$ is a Pochhammer symbol and $n!=1\cdot2\cdot...\cdot n$.

It is well-known that the series $_{p}F_{q}$ converges for all $x$ if $p\leq q$ and for $|x|<1$ if $p=q+1$, and the series diverges 
for all $x\neq 0$ if $p>q+1$. The convergence of the series for the case $|x|=1$ when $p=q+1$ is of much interest.
For the convergence of the series in this case we have the following theorem (Andrews, G., Askey, R., \& Roy, R., 2000)
\begin{theorem} \label{hypergeom}
The series $_{q+1}F_{q}=(a_{1},...,a_{q+1}; b_{1},...,b_{q};x)$ with $|x|=1$ converges absolutely if $\left(\sum b_{i} - \sum a_{i}\right) >0$.
The series converges conditionally if $x=e^{i\theta}\neq 1$ and $0\geq \left(\sum b_{i} - \sum a_{i}\right)>-1$. 
The series diverges if $\left(\sum b_{i} - \sum a_{i}\right) \leq -1$.
\end{theorem}

And maybe one of the most significant hypergeometric series that have numerous applications is the $_{2}F_{1}(\alpha,\beta;\gamma;x)$ hypergeometric functions, 
Gauss' hypergeometric series
\begin{eqnarray}
_{2}F_{1}(\alpha,\beta;\gamma;x)  \label{hyper} \nonumber
&=&
\sum_{n=0}^{\infty}\frac{(\alpha)_{n}(\beta)_{n}}{n!(\gamma)_{n}}x^{n}= \\ 
&=&
\sum_{n=0}^{\infty}\frac{\alpha\cdot(\alpha+1)\cdot...\cdot(\alpha+n-1)\cdot\beta\cdot(\beta+1)\cdot...
\cdot(\beta+n-1)}{n!\cdot\gamma\cdot(\gamma+1)\cdot...\cdot(\gamma+n-1)}x^{n} = \\ \nonumber
&=&
1+\frac{\alpha\cdot\beta}{1\cdot\gamma}x+
\frac{\alpha\cdot(\alpha+1)\cdot\beta\cdot(\beta+1)}{1\cdot2\cdot\gamma\cdot(\gamma+1)}x^{2}+
\frac{\alpha\cdot(\alpha+1)\cdot(\alpha+2)\cdot\beta\cdot(\beta+1)\cdot(\beta+2)}{1\cdot2\cdot3\cdot\gamma\cdot(\gamma+1)\cdot(\gamma+2)}x^{3}+... \nonumber
\end{eqnarray}
defined for $|x|<1$ and by analytic continuation elsewhere. 
We see from theorem \ref{hypergeom} that the series $_{2}F_{1}(\alpha,\beta;\gamma;x)$ diverges in general for $x=1$ and $\gamma-\alpha-\beta\leq 0$ and for $x=-1$ and $\gamma-\alpha-\beta\leq -1$.

One of the most important summation formulas for $_{2}F_{1}$ hypergeometric series is given by the binomial theorem:
$$
_{2}F_{1}(\alpha,\beta;\beta;x)=\:_{1}F_{0}(\alpha;-;x)=\sum_{n=0}^{\infty}\frac{(\alpha)_{n}}{n!}x^{n}=(1-x)^{-\alpha},
$$
where $|x|<1$.

Since the binomial theorem is at the foundation of most of the summation formulas for hypergeometric series (Gasper, G. \& Rahman, M., 2004), 
we shall conduct our reasoning in relation to the binomial series as well.
The binomial series is one of the most significant and important infinite series which often occurs in applications and in different areas of mathematics
\begin{equation}
1+\alpha x+\frac{\alpha(\alpha-1)}{1\cdot 2}x^{2}+\frac{\alpha(\alpha-1)(\alpha-2)}{1\cdot 2\cdot 3}x^{3}+...+\frac{\alpha(\alpha-1)(\alpha-2)...(\alpha-n+1)}{1\cdot 2\cdot 3...\cdot n}x^{n}+... \label{series}
\end{equation}
which is the Taylor series for the function $f(x)=(1+x)^{\alpha}$, where $\alpha$ is a real number. When $\alpha$ is a natural we get the expansion having finite
number of terms, a well-known Newton binomial formula. 

Although the binomial series always converges when $|x|<1$, the question of whether
or not it converges at the endpoints $x=\pm 1$, depends on the value of $k$ and should be checked every time for concrete series under study. 
For this reason, the applicability of the binomial series (\ref{series}) for the values of $x=\pm 1$ in a general case is always under question, 
and requires additional analysis of a series.

However, the binomial series (\ref{series}) can be expressed by means of the hypergeometric function $_{2}F_{1}(\alpha,\beta;\gamma;x)$
if putting $\alpha=-m$, $\beta=\gamma$, and replacing $x$ by $-x$. 
So, from the behavior of the hypergeometric series (\ref{hyper}) at the endpoints $x=\pm 1$, it turns out that the binomial series converges at $x=1$ if $-1<m<0$ and $m>0$, 
and at $x=-1$ if $m>0$. We see that the condition $-1<m<0$ is still a restriction for the convergence of the series at $x=1$.

In this short note, within a new theoretical framework, we establish,  without any use of analytic continuation technique, that the hypergeometric series $_{2}F_{1}(-\alpha,\beta;\beta;-x)$, that is the binomial series (\ref{series}), is also valid for $x=1$ 
for any integer power $\alpha$ not equal to zero, without imposing any restrictions on $\alpha$.

\section{Main Results}\par

We have
\begin{theorem} \label{mytheorem}
For any integer $a\neq 0$ and real $x$, $-1<x\leq 1$, the following equality holds
\begin{eqnarray}
_{2}F_{1}(-a,\beta;\beta;-x) &=& (1+x)^a =
\sum_{u=0}^{\infty}x^{u}\:\biggl(\:\begin{matrix}a\\ u\end{matrix}\:\biggr)= 1+ax+\frac{a(a-1)}{1\cdot 2}x^{2}+...+  \\ \nonumber
&+&
\frac{a(a-1)(a-2)...(a-n+1)}{1\cdot 2\cdot 3...\cdot n}x^{n}+... \label{myformula}
\end{eqnarray}
where
$$
\left(\:\begin{matrix}a\\ u\end{matrix}\:\right)=\frac{a(a-1)(a-2)...(a-u+1)}{1\cdot 2\cdot 3...\cdot u}, \;\;\; (u\geq 1), \;\;\; 
\left(\:\begin{matrix}a\\ 0\end{matrix}\:\right)=1
$$
\end{theorem}

For preparation of proof, we need several auxiliary propositions. So, for the sake of logical completeness of our reasoning, we recall here some facts from 
(Varshamov, R. \& Bagdasaryan, A., 2009); we do not give the proofs of all propositions and refer the interested reader to our work mentioned above.
The proofs will be given only to the theorems that we directly use to prove our main result of theorem \ref{mytheorem}.

In (Varshamov, R. \& Bagdasaryan, A., 2009) we introduced a new method for integers ordering, 
so that the integer numbers are reorganized in such a way that negative numbers become beyond infinity, 
that is positive and negative numbers are linked through infinity, and the set of integers $\mathbb{Z}=[0,1,2...,-2-1]$. 
Recently, based on this method of ordering, we presented a new and elementary approach to values of the 
Riemann zeta function at non-positive integer points (Bagdasaryan, A., 2008).
So, let us introduce necessary definitions and statements of our new theoretical setting that we use in this work.

\begin{definition} 
We shall say that $a$ precedes $b$, $a, b \in \mathbb{Z}$, and write $a\prec b$, if and only if $\frac{-1}{a}<\frac{-1}{b}$; $a\prec b \Leftrightarrow \frac{-1}{a}<\frac{-1}{b}$ \footnote{by convention $0^{-1}=\infty$}.
\end{definition}

For the newly introduced ordering of integers, the following two essential axioms of order are true:
\begin{itemize}
	\item Transitivity \newline
					$\forall a, b$: if $a\prec b$ and $b\prec c$ then $a\prec c$
	\item Connectedness \newline
					$\forall a, b$: if $a\neq b$ then $a\prec b$ or $b\prec a$
\end{itemize}

Thus, this new ordering on $\mathbb{Z}$ is a strict total (linear) order.

We also introduced a new axiomatic system with arithmetic as a model, and defined a new class of regular functions.

\begin{definition} \label{regdef}
The function $f(x)$, $x\in \mathbb{Z}$ , is called \textit{regular} if there exists an elementary \footnote{determined by formulas containing a finite number of algebraic or trigonometric operations performed over argument, function and constants} function $F(x)$ such that $F(z+1)-F(z)=f(z), \; \forall z\in \mathbb{Z}$. The function $F(x)$ is said to be a \textit{generating function} for $f(x)$.
\end{definition}

From Definition \ref{regdef} of regular functions the following equality is deduced

\begin{equation}
F(b+1)-F(a)=\sum_{u=a}^{b}f(u), \label{regeq}
\end{equation}

which is true for any values of $a$ and $b$, $a\leq b$.

\textbf{Remark}. One can note an analogy between the equality (\ref{regeq}) and the problem of summation of functions 
in the theory of finite differences. Roughly speaking, the problem of finding the sum $S_{n}=f(0)+f(1)+...+f(n)$
as a function of $n$ can be reduced to finding the function $F(x)$ such that $F(x+1)-F(x)=f(x)$ for a given function $f(x)$.
Indeed, if such a function $F(x)$ is found then for $x$ taking consecutively values $0,1,...,n$ one has $F(1)-F(0)=f(0)$, $F(2)-F(1)=f(1)$,...,
$F(n+1)-F(n)=f(n)$. Adding these equalities we arrive at $f(0)+f(1)+...+f(n)=\sum_{j=0}^{n}f(j)=F(n+1)-F(0)$, which is the solution of the problem of
summation of function $f(x)$ (Gel'fond, A. O., 1967). 
Taking as the limits of summation $b$ instead of $n$ and $a$ instead of $0$, we come to (\ref{regeq}). 

Let $f(x)$ be a function of real variable defined on $\mathbb{Z}$.  

\begin{definition}
For any integer numbers $a, b\in \mathbb{Z}$
$$
\sum_{u=a}^b{f(u)}=\sum_{u\:\in\: \mathbb{Z}_{a, b}}{f(u)}  
$$
where $\mathbb{Z}_{a, b}$ is a part of $\mathbb{Z}$ such that
\begin{equation} \label{def:sum}
\mathbb{Z}_{a, b} = \left\{              
									\begin{array}{ll}                   
									[a, b] & (a\preceq b)\\                   
									\mathbb{Z}\setminus (b, a) & (a\succ b)              
									\end{array}       
			 				\right.
\end{equation}
and $\mathbb{Z}\setminus (b, a)=[a, -1]\cup[0, b]$ 
\end{definition}

The Definition 5 extends the classical definition of sum to the case when $b<a$.

The set $\mathbb{Z}_{a, b}$, depending on the elements $a$ and $b$, can be both finite and infinite. Thus, the sum on the right-hand side of (\ref{def:sum}) may become an infinite series.  
Following L. Euler\footnote{who was convinced that ``to every series could be assigned a number'', which is the first Euler's principle on infinite series}, 
we take as a postulate the assertion that \textit{any series $\sum_{u=1}^{\infty}f(u)$,
where $f(u)$ is a regular function, has a certain finite, numeric, value}.

The regular functions are required to satisfy the following quite natural system of "axioms".
\begin{enumerate}
\item If $S_n=\sum_{u=a}^n{f(u)} \;\; \forall n$, then $\lim_{n\rightarrow \infty}S_n=\sum_{u=a}^\infty {f(u)}$, \; ($n\rightarrow \infty$ means that $n$ unboundedly increases, without changing the sign).\\
\item If $S_n=\sum_{u=1}^{n/2}{f(u)} \;\; \forall n$, then $\lim_{n\rightarrow \infty}S_n=\sum_{u=1}^\infty {f(u)}$.\\
\item If $\sum_{u=a}^\infty {f(u)}=S$, then $\sum_{u=a}^\infty {\lambda f(u)}=\lambda S, \; \lambda\in R$.\\
\item  If $\sum_{u=a}^\infty {f_1(u)}=S_1$ and $\sum_{u=a}^\infty {f_2(u)}=S_2$, then  $\sum_{u=a}^\infty {\left(f_1(u)+f_2(u)\right)}= S_1+S_2$.\\
\item For any $a$ and $b$, $a\leq b$: $F(b+1)-F(a)=\sum_{u=a}^b{f(u)}$.\\
\item If $G=[a_1, b_1]\cup[a_2, b_2]$, $[a_1, b_1]\cap[a_2, b_2]=\emptyset$, then $\sum_{u\in G}{f(u)}=\sum_{u=a_1}^{b_1}{f(u)}+\sum_{u=a_2}^{b_2}{f(u)}$. \\
\end{enumerate}
 
It is worth to notice that introduced ``axioms'' are consistent\footnote{that is they do not contradict to
what has been obtained in the framework or on the basis of known means of analysis}
with statements and definitions of analysis: the axiom 1 is the classical definition of sum of a series; the axiom 2 also does not contradict to classical analysis; the axioms 3 and 4 constitute the well-known properties of convergent series\footnote{which have a finite sum; bearing in mind our postulate, usage of these axioms are in agreement with classical analysis}; the axiom 5 was mentioned above; the axiom 6 has quite common formal meaning.

This axiomatic system defines the method for summation of infinite series, from which we obtain a general and unified approach to summation of divergent series.
The method is regular, since, due to axiom 1, every convergent series is summable to its usual sum.

\begin{definition} \label{def:limit_seq}
A number $A$ is said to be the \textit{limit} of the numeric sequence $F(1)$, $F(2)$,..., $F(n)$,... (function of integer argument), i. e. $\lim_{n\rightarrow\infty}F(n)=A$, if $\sum_{u=1}^{\infty}f(u)=A$, where $F(1)=f(1)$ and $F(u)-F(u-1)=f(u) \; \; (u>1)$, that is
$$
F(1)+(F(2)-F(1))+...+\left(F(n)-F(n-1)\right)+...=A. 
$$
\end{definition}

\textbf{Remark}.
The Definition \ref{def:limit_seq} in particular coincides with the classical definition of limit of a sequence, which is convergent in usual sense, and reduces the problem of existence of the limit of functions of integer argument $F(n)$ to finding the sum of the series
$$
F(1)+(F(2)-F(1))+...+(F(n)-F(n-1))+...
$$
In force of our postulate, any elementary function of integer argument defined on $\mathbb{Z}$ has a certain limit. 
This, in turn, allows us to extend the classical theory of limits and determine the limits of unbounded and oscillating sequences and functions.

In our proof of theorem \ref{mytheorem} we consecutively rely on the propositions below, 
which we obtain from the introduced axiomatic system. Some of the propositions, mainly those concerning limits of unbounded and oscillating functions and sequences, 
which are true with the new number line and within the framework of our definitions, may be false in the usual (``epsilon-delta'')-definition of limit in classical analysis, 
and also since as it is known real and complex analysis fail in finding the limits of sequences/functions with oscillation or unboundedness;
the matter is that the topology\footnote{but not axiomatics, which is obviously consistent with the standard topology of $\mathbb{R}$} 
of the newly ordered number line may be different\footnote{which is the subject of a separate study} from that in classical sense. Probably, section 2 and
lemmas 2.15, 2.16 and their corollaries in (Varshamov, R. \& Bagdasaryan, A., 2009) may shed some light on the topological properties of the new number line.

\begin{proposition}
Let $\mu(x)$ be any even elementary function, $t$ is a fixed natural number, $\epsilon=\pm 1$, and $\delta=(1-\epsilon)/2$. Then
\begin{equation}
\lim_{n\rightarrow\infty}(-1)^{n}\sum_{u=\delta}^{t-1+\delta}\mu\left(n+\frac{\epsilon t}{2}-\epsilon u\right)=0   \label{corollary}
\end{equation}
\end{proposition}

\begin{theorem} \label{theorsum}
For any regular function $f(x)$ such that $f(-x)=f(x-\epsilon t)$, $\epsilon=\{0, \pm 1\}$, we have
$$
\sum_{u=1}^{\infty}f(u)=\frac{\epsilon}{2} \sum_{u=\delta}^{t-1+\delta}\left(\lim_{n\rightarrow\infty}f(n-\epsilon u)-f(-\epsilon u)\right)-\frac{1}{2}f(0)
$$
\end{theorem}

\begin{proposition}
For any natural number $k$
\begin{equation}
\lim_{n\rightarrow\infty}(-1)^{n}(2n+1)^{k}=0 \label{limzero}
\end{equation}
\end{proposition}
\textbf{Proof}. Let us consider the cases of $k$ even and $k$ odd.
\begin{enumerate}
	\item Let $k$ is even, $k=2m$. \newline  
	From the formula (\ref{corollary}), for $\epsilon=1$ and $t=1$, we have
	\begin{equation}
	\lim_{n\rightarrow\infty}(-1)^{n}\mu\left(n+\frac{1}{2}\right)=0. 
	\end{equation}
	Taking $\mu(x)=\left(2x\right)^{2m}$, we get
	$$
	\lim_{n\rightarrow\infty}(-1)^{n}\mu\left(n+\frac{1}{2}\right)=\lim_{n\rightarrow\infty}(-1)^{n}\left(2n+1\right)^{2m}=0
	$$
	and
	$$
		\lim_{n\rightarrow\infty}(-1)^{n}(2n+1)^{k}=0.
	$$
	\item Let $k$ is odd, $k=2m-1$. \newline
	From Theorem \ref{theorsum} for $\epsilon=-1$ and $t=1$, we get
	\begin{equation}
	\sum_{u=1}^{\infty}f(u)=-\frac{1}{2}\lim_{n\rightarrow\infty}f(n+1)   \label{theorcor}
	\end{equation}
	The function $f(x)=(-1)^{x}\left(2x-1\right)^{2m-1}$ satisfies the condition of Theorem \ref{theorsum}. Then, substituting $(-1)^{x}\left(2x-1\right)^{2m-1}$ into (\ref{theorcor}), we obtain
	$$
	\sum_{u=1}^{\infty}(-1)^{u-1}\left(2u-1\right)^{2m-1}=-\frac{1}{2}\lim_{n\rightarrow\infty}(-1)^{n}\left(2n+1\right)^{2m-1}.
	$$
	From the other hand, from Theorem \ref{theorsum} for $\epsilon=0$, that is for even regular functions, we have
	$$
	\sum_{u=1}^{\infty}f(u)=-\frac{1}{2}f(0).
	$$
	Hence, taking as a generating function, the function 
	$$
	F(n)=(-1)^{n}\left(\sum_{u=1}^{m}\beta_{u}\left(n-\frac{1}{2}\right)^{2u-1}\sin\left(n-\frac{1}{2}\right)\theta +
				\sum_{u=0}^{m-1}\overline{\beta}_{u}\left(n-\frac{1}{2}\right)^{2u}\cos\left(n-\frac{1}{2}\right)\theta\right),
	$$
	where the coefficients $\beta_{u}$ and $\overline{\beta}_{u}$ are taken in such a way, so that the equality
	$F(n+1)-F(n)=(-1)^{n-1}n^{2m-1}\sin n\theta$ holds, we get 
	$$
	\sum_{u=1}^{\infty}(-1)^{u-1}u^{2m-1}\sin u\theta=0, \;\;\;\;\; -\pi<\theta<\pi \;\;\;\;\;\;\;\;\;\;\;\;\;\;\;\;\;\;\;\;\;\;\; \textup{(Hardy, G. H., 1949)}
	$$
	Then, putting $\theta=\frac{\pi}{2}$, we arrive at
	$$
	1^{2m-1}-3^{2m-1}+5^{2m-1}-7^{2m-1}+...=\sum_{u=1}^{\infty}(-1)^{u-1}(2u-1)^{2m-1}=0.
	$$
	Therefore, we obtain
	$$
	\sum_{u=1}^{\infty}(-1)^{u-1}\left(2u-1\right)^{2m-1}=-\frac{1}{2}\lim_{n\rightarrow\infty}(-1)^{n}\left(2n+1\right)^{2m-1}=0
	$$
	and
	$$
	\lim_{n\rightarrow\infty}(-1)^{n}\left(2n+1\right)^{k}=0.
	$$
\end{enumerate}
The theorem is proved completely.

\begin{theorem} \label{polynom}
Suppose $f(x)$ is a polynomial defined over the field of real numbers, $x\in\mathbb{R}$. Then
$$
\lim_{n\rightarrow\infty}(-1)^{n}f(n)=0
$$
\end{theorem}
\textbf{Proof}. Because the limit of algebraic sum of finite number of sequences equals to
the algebraic sum of limits of the sequences, which is derived from axiom 1, that is
$$
\lim_{n\rightarrow\infty}\sum_{u=1}^{m}\alpha_{u}F_{u}(n)=\sum_{u=1}^{m}\alpha_{u}\lim_{n\rightarrow\infty}F_{u}(n),
$$
where $\alpha_{u}$ are real numbers, we have
\begin{eqnarray}
\lim_{n\rightarrow\infty}(-1)^{n}f(n)&=&\lim_{n\rightarrow\infty}(-1)^{n}\left(a_{k}n^{k}+a_{k-1}n^{k-1}+...+a_{1}n+a_{0}\right)= \nonumber \\
&=& 
a_{k}\lim_{n\rightarrow\infty}(-1)^{n}n^{k}+a_{k-1}\lim_{n\rightarrow\infty}(-1)^{n}n^{k-1}+...+a_{1}\lim_{n\rightarrow\infty}(-1)^{n}n+
a_{0}\lim_{n\rightarrow\infty}(-1)^{n}. \nonumber
\end{eqnarray}

To prove the theorem, it is sufficient to show that for every non-negative integer value $\sigma$
\begin{equation}
\lim_{n\rightarrow\infty}(-1)^{n}n^{\sigma}=0. \label{sigma}
\end{equation}
The proof is by induction over $\sigma$. 
\begin{enumerate}
	\item Let $\sigma=0$. From the formula (\ref{corollary}), as a special case for $\epsilon=1$ and $t=1$, it follows that
	\begin{equation}
	\lim_{n\rightarrow\infty}(-1)^{n}\mu\left(n+\frac{1}{2}\right)=0. \label{sigmanull}
	\end{equation}
Then, according to (\ref{sigmanull}), putting $\mu\left(n+\frac{1}{2}\right) \equiv 1$, we immediately obtain that formula (\ref{sigma}) holds true for $\sigma=0$.
	\item Assume now that (\ref{sigma}) holds for all positive integers less than some natural number $k$, $\sigma<k$. 
	Then, from (\ref{limzero}), applying the binomial theorem to $(2n+1)^{k}$, we get the expansion
	$$
	2^{k}\lim_{n\rightarrow\infty}(-1)^{n}n^{k}+2^{k-1}\biggl(\:\begin{matrix}k\\ 1\end{matrix}\:\biggr)\lim_{n\rightarrow\infty}(-1)^{n}n^{k-1}+
	2^{k-2}\biggl(\:\begin{matrix}k\\ 2\end{matrix}\:\biggr)\lim_{n\rightarrow\infty}(-1)^{n}n^{k-2}+...+\lim_{n\rightarrow\infty}(-1)^{n}=0,
	$$
	in which, by the inductive assumption, all terms, except for the first, are equal to zero. Therefore, the first term will be equal to zero as well
	$$
	2^{k}\lim_{n\rightarrow\infty}(-1)^{n}n^{k}=0
	$$
	and finally
	$$
	\lim_{n\rightarrow\infty}(-1)^{n}n^{k}=0
	$$
	and the formula (\ref{sigma}) is also true for $\sigma=k$.
	
	Thus, by virtue of the induction, we obtain that the formula (\ref{sigma}) holds for all non-negative integer values $\sigma$. The theorem is proved.
\end{enumerate}

\textbf{Proof of Theorem \ref{mytheorem}}. For $a>0$ the $_{2}F_{1}$ becomes a finite sum, that is the formula (\ref{myformula}) is reduced to the elementary binomial formula. Thus, let $a$ be a negative number, that is $a=-m$. In this case, dividing $1$ by $(1+x)^m$, using the rule of division of one polynomial by another, we will gradually obtain in the quotient the terms of the right-hand side of (\ref{myformula}). Stopping the process of division at some fixed number $k$, we get

$$
(1+x)^{-m} = 1-x\biggl(\:\begin{matrix}m\\ 1\end{matrix}\:\biggr)+x^{2}\biggl(\:\begin{matrix}m-1\\ 2\end{matrix}\:\biggr)-...(-1)^{k}x^{k}\biggl(\:\begin{matrix}m-k-1\\ k\end{matrix}\:\biggr)+R_{k}^{m}(x),
$$

where

$$
R_{k}^{m}(x)=\frac{(-1)^{k+1}}{(1+x)^m}\sum_{u=0}^{m-1}\biggl(\:\begin{matrix}k+u\\ u\end{matrix}\:\biggr)\biggl(\:\begin{matrix}m+k\\ m-1-u\end{matrix}\:\biggr)x^{k+1+u}.
$$

If $|x|<1$ and $k$ is large enough, then obviously $R_{k}^{m}(x)$ is arbitrarily small and it tends to zero as $k$ unboundedly increases, that is

$$
\lim_{k\rightarrow\infty}R_{k}^{m}(x)=0 \;\;\;\;\;\;\; |x|<1.
$$

If $x=1$, then $\bigl|R_{k}^{m}(x)\bigr|$ becomes a polynomial of degree $m-1$ in the variable $k$, which according to the statement of theorem \ref{polynom} also satisfy the equality

$$
\lim_{k\rightarrow\infty}R_{k}^{m}(1)=0.
$$

This means that (\ref{myformula}) holds for all $x$, $-1<x\leq 1$.\\
The proof is completed.

\textbf{Example}. Here we present several examples to our theorem.
\begin{eqnarray}
	_{2}F_{1}(1,\beta;\beta;-1) &=& (1+1)^{-1}=1-1+1-...=\frac{1}{2}  \nonumber \;\;\;\;\;\;\;\;\;\;\;\;\;\;\;\;\;\;\;\;\;\;\;\;\;\; \textup{(Hardy, G. H., 1949)} \\
	_{2}F_{1}(2,\beta;\beta;-1) &=& (1+1)^{-2}=1-2+3-...=\frac{1}{4} \nonumber \;\;\;\;\;\;\;\;\;\;\;\;\;\;\;\;\;\;\;\;\;\;\;\;\;\; \textup{(Hardy, G. H., 1949)} \\
	_{2}F_{1}(3,\beta;\beta;-1) &=& (1+1)^{-3}=1-3+6-...=\frac{1}{8}  \nonumber \;\;\;\;\;\;\;\;\;\;\;\;\;\;\;\;\;\;\;\;\;\;\;\;\;\; \textup{(Hardy, G. H., 1949)} 
\end{eqnarray}
and so on. 

\textbf{Remark}. The validity of the Theorem \ref{mytheorem} can also be established with use of the formula 
$\left(\:\begin{matrix}-n\\ m\end{matrix}\:\right)=(-1)^{m}\left(\:\begin{matrix}n+m-1\\ m\end{matrix}\:\right)$ (Riordan, J., 1958)
and Theorem \ref{polynom}.

\section{Conclusion}\par

In this paper we used a new theoretical setting to reprove the convergence of the hypergeometric series
$_{2}F_{1}$ for the special values of parameters and for $|x|<1$, and to prove its convergence at the special argument of $-1$. 

The general setting that we presented here is a useful analytical tool
as it enables one to derive, prove, and treat known mathematical facts and obtain new results without applying analytic continuation techniques.  
Another  advantage is that it facilitates to get simpler and elementary solutions to problems and to find qualitatively new results.

In particular, we plan to continue studying the special cases of the hypergeometric function $_{2}F_{1}$ and in general $_{p}F_{q}$ for some other special values of $p$ and $q$ within the framework of this new theoretical setting in our future works.

\textbf{Acknowledgment} \par
It is a pleasure to thank the anonymous referees for valuable comments and suggestions leading to an improvement of the
presentation of this paper. 


\textbf{References} \par

\par
Andrews, G., Askey, R., and Roy, R. (2000). \emph{Special functions}. Cambridge: Cambridge University Press.\par
Gasper, G. \& Rahman, M. (2004). \emph{Basic hypergeometric series}. Cambridge: Cambridge University Press. \par
Gel'fond, A. O. (1967). Calculus of Finite Differences. 
Moscow: Nauka, 3rd rev. ed.; 
German transl., 1958. Berlin: VEB Deutscher Verlag. (Transl. of 1st Russian Ed.); 
French transl., 1963. Paris: Collection Univ. de Math., XII, Dunod. (Transl. of 1st Russian Ed.); 
English transl., 1971. Delhi: Hindustan Publ. (Transl. of 3rd rev. Russian Ed.). \par
Hardy, G. H. (1949). \emph{Divergent series}. Oxford: Oxford University Press. \par
Riordan, J. (1958). \emph{An introduction to combinatorial analysis}. New York: John Wiley \& Sons. \par
Temme, N. M. (1996), \emph{Special functions: an introduction to the classical functions of mathematical physics}. New York: John Wiley \& Sons. \par
Varshamov, R. \& Bagdasaryan, A. (2009), On one number-theoretic conception: towards a new theory, submitted. [Online] Available: arXiv:0907.1090 [math.GM],  \href{http://arxiv.org/abs/0907.1090}{http://arxiv.org/abs/0907.1090} (July 6, 2009) \par
Bagdasaryan, A. (2008), An elementary and real approach to values of the Riemann zeta function. Proceedings of the XXVII International
Colloquium on Group Theoretical Methods in Physics, Yerevan, Armenia, August 13-19, 2008. \emph{Physics of Atomic Nuclei}, (to appear). 
[Online] Available: arXiv:0812.1878 [math.NT], \href{http://arxiv.org/abs/1812.1878}{http://arxiv.org/abs/0812.1878} (v1: December 10, 2008; v2: July 23, 2009) \par
Whittaker, E. T. \& Watson, G. N. (2005), \emph{A course of modern analysis}. Cambridge: Cambridge University Press. \par
\bigskip




\newpage
\end{document}